\documentclass[11pt]{amsart}

\usepackage{amsmath,amsfonts,amsthm,amssymb, mathtools}
\usepackage{setspace}
\usepackage{ mathrsfs }
\usepackage{ marvosym }
\usepackage{stackengine}
\usepackage{graphicx}
\usepackage{bm}
\usepackage{caption}
\usepackage{subfig}
\usepackage{floatrow}

\usepackage{Tabbing}
\usepackage{fancyhdr}
\usepackage{lastpage}
\usepackage{extramarks}

\usepackage{chngpage}
\usepackage{soul,color}
\usepackage{graphicx,float,wrapfig}
\usepackage{fullpage}
\usepackage{hyperref}
\usepackage[all]{xy}

\usepackage{enumerate}
\usepackage{marvosym}
\usepackage{enumitem}
\newcommand{\be}{\begin{enumerate}}
\newcommand{\ee}{\end{enumerate}}
\newcommand{\beq}{\begin{equation}}
\newcommand{\eeq}{\end{equation}}
\usepackage{verbatim}

\newcommand{\ba}{\begin{align*}}

\newcommand{\ea}{\end{align*}}

\newcommand{\R}{{\mathbb R}}

\newcommand{\T}{{\mathcal{T}}}

\newcommand{\Z}{{\mathbb Z}}

\newcommand{\C}{{\mathbb C}}

\renewcommand{\H}{{\mathbb H}}

\newcommand{\GL}{\operatorname{GL}}

\newcommand{\SL}{\operatorname{SL}}

\newcommand{\g}{\mathbf{g}}

\renewcommand{\t}{\mathfrak{t}}

\newcommand{\im}{\text{Im}}

\newcommand{\Hom}{\text{Hom}}

\newcommand{\id}{\text{id}}
\newcommand{\ol}[1]{\overline{#1}}

\newcommand{\pd}[2]{\frac{\partial#1}{\partial#2}}

\newcommand{\ra}{\rightarrow}
\newcommand{\tra}\twoheadrightarrow
\newcommand{\hra}{\hookrightarrow}

\newcommand{\tla}\twoheadleftarrow

\newcommand{\vn}{\varepsilon}

\newcommand{\lp}{\left(}

\newcommand{\rp}{\right)}
\newcommand{\lpi}{\left|}
\newcommand{\rpi}{\right|}

\newcommand{\lbrac}{\left[}

\newcommand{\rbrac}{\right]}

\newcommand{\llambda}{\bm\lambda}
\newcommand{\mmu}{\bm\mu}
\renewcommand{\v}{\mathbf{v}}
\renewcommand{\t}{\mathbf{t}}
\renewcommand{\Re}{\text{Re}}
\renewcommand{\Im}{\text{Im}}
\newcommand{\E}{\mathcal{E}}
\newcommand{\MCG}{\text{MCG}}

\def\bea#1\eea{\begin{align*}#1\end{align*}}
\def\bc#1\ec{\begin{comment}#1\end{comment}}

\newtheorem{Theorem}{Theorem}[section]

\newtheorem{Lemma}[Theorem]{Lemma}
\newtheorem{Proposition}[Theorem]{Proposition}
\newtheorem{Corollary}[Theorem]{Corollary}
\newtheorem{Conjecture}[Theorem]{Conjecture}

\newcommand{\Implies}[2]{$\text{\ref{#1}}\implies\text{\ref{#2}}$}

\let\inf\relax \DeclareMathOperator*\inf{\vphantom{p}inf}
\let\lim\relax \DeclareMathOperator*\lim{lim\vphantom{p}}

\title{Classifying complex geodesics for the Carath\'eodory metric on low-dimensional Teichm\"uller spaces}
\date{}
\author{Dmitri Gekhtman and Vladimir Markovic}

\begin{document}
\captionsetup[subfigure]{labelfont=rm} 
\captionsetup[subfloat]{captionskip = 5pt}
\captionsetup{font=footnotesize}

\maketitle

\begin{abstract}
It was recently shown that the Carath\'eodory and Teichm\"uller metrics on the Teichm\"uller space of a closed surface do not coincide.
On the other hand, Kra earlier showed that the metrics coincide when restricted to a Teichm\"uller disk generated by a differential with no odd-order zeros.
Our aim is to classify Teichm\"uller disks on which the two metrics agree, and we conjecture that the Carath\'eodory and Teichm\"uller metrics agree on a Teichm\"uller disk if and only if the Teichm\"uller disk is generated by a differential with no odd-order zeros.
Using dynamical results of Minsky, Smillie, and Weiss, we show that it suffices to consider disks generated by Jenkins-Strebel differentials. We then prove a complex-analytic criterion characterizing Jenkins-Strebel differentials which generate disks on which the metrics coincide. Finally, we use this criterion to prove the conjecture for the Teichm\"uller spaces of the five-times punctured sphere and the twice-punctured torus.
We also extend the result that the Carath\'eodory and Teichm\"uller metrics are different to the case of compact surfaces with punctures.  
\end{abstract}

\section{Introduction}
Let $\T:=\T_{g,n}$ denote the Teichm\"uller space of a finite-type orientable surface $S_{g,n}$. Let $\H$
denote the upper half-plane, equipped with its Poincar\'e metric $d_\H$.
The Carath\'eodory metric on $\T$ is the smallest metric so that every holomorphic map $\T \ra (\H, d_\H)$ is nonexpanding.
On the other hand, the Kobayashi metric on $\T$ is the largest metric so that every map $(\H,d_\H) \ra \T$ is nonexpanding.
Royden \cite{Ro} proved that the Kobayashi metric is the same as the classical Teichm\"uller metric.
Whether or not the Carath\'eodory metric is also the same as the Teichm\"uller metric was a longstanding open problem.

Let $\tau:\H \ra \T$ be a Teichm\"uller disk.
Then the Carath\'eodory and Kobayashi metrics agree on $\tau(\H)$ if and only if there is a {\em holomorphic retraction}
onto $\tau(\H)$, i.e. a holomorphic map $F:\T \ra \H$ so that $F\circ \tau = \id_\H$. Thus, the problem of determining whether
the Carath\'eodory and Kobayashi metrics agree reduces to checking whether each Teichm\"uller disk is a holomorphic retract
of Teichm\"uller space. In 1981, Kra \cite{Kr} showed that if a holomorphic quadratic differential has no odd-order zeros, then its associated
Teichm\"uller disk is a holomorphic retract. However, it was recently shown \cite{Mar} that not all Teichm\"uller disks in $\T_g$ are retracts,
and so the Carath\'eodory and Kobayashi metrics are different. 

It remains to classify the Teichm\"uller disks on which the two metrics agree. In other words, we would like to know which Teichm\"uller disks are holomorphic retracts of Teichm\"uller space. Put another way, our aim is to classify the complex geodesics for the Carath\'eodory metric on Teichm\"uller space.
 
We conjecture the converse of Kra's result:

\begin{Conjecture}\label{Conj}
A Teichm\"uller disk is a holomorphic retract if and only if it is generated by a quadratic differential with no odd-order zeros.
\end{Conjecture}

In this paper, we suggest a program for proving the conjecture. We carry out the program for the spaces $\T_{0,5}$ and $\T_{1,2}$ of complex dimension two. That is, we prove
\begin{Theorem}[\bf{Main Result}]\label{Main}
A Teichm\"uller disk in $\T_{0,5}$ or $\T_{1,2}$  is a holomorphic retract if and only if the zeros of $\phi$
are all even-order.
\end{Theorem}

Dynamics on the moduli space of quadratic differentials plays a key role in the proof.
If the quadratic differential $\phi$ generates a holomorphic retract, then so does every differential in its $\SL_2(\R)$ orbit closure.
On the other hand, combined results of Minsky-Smillie \cite{Mi} and Smillie-Weiss \cite{Sm} show that
the orbit closure contains a Jenkins-Strebel differential in the same stratum as $\phi$. To prove Conjecture \ref{Conj},
it thus suffices to consider Jenkins-Strebel differentials. To this end, we develop a complex-analytic criterion, Theorem \ref{criterion}, characterizing Jenkins-Strebel differentials which generate retracts. The criterion involves a certain holomorphic embedding $\E^\phi:\H^k \ra \T$, called
a {\em Teichm\"uller polyplane}, of the $k$-fold product of $\H$ into Teichm\"uller space. 
To prove Conjecture $\ref{Conj}$ for $\T_{0,5}$, we reduce to the case that $\phi$ is an {\em L-shaped pillowcase}. Then, following the argument in \cite{Mar}, we use our analytic criterion to show that an $L$-shaped pillowcase does not generate a retract.            
   
\subsection{The Carath\'eodory and Kobayashi pseudometrics}
 
 A {\em Schwarz-Pick system} is a functor assigning to each complex manifold $M$ a pseudometric $d_M$ satisfying the following conditions:
\begin{enumerate}[label=(\roman*)]
\item
The metric assigned to the upper half-plane $\H = \{\lambda \in \C | \im \lambda>0 \}$ is the Poincar\'e metric of curvature $-4$:
$$
d_\H(\lambda_1,\lambda_2) = \tanh^{-1} \lpi \frac{ \lambda_1-\lambda_2 }{ \lambda_1-\ol{\lambda}_2}\rpi .
$$
\item
Any holomorphic map $f:M\ra N$ between complex manifolds is non-expanding:
$$
d_N\lp f(p), f(q)\rp \leq d_M(p,q)
$$ 
for all $p,q \in M$.
\end{enumerate}

Distinguished among the Schwarz-Pick systems are the Carath\'eodory and Kobayashi pseudometrics.
The Carath\'eodory pseudometric $C_M$ on a complex manifold $M$ is the smallest pseudometric so that all holomorphic maps
from $M$ to $\H$ are nonexpanding.
More explicitly, 
\begin{equation}\label{D1}
C_M(p,q) = \sup_{f \in \mathcal{O}(M, \H)} d_\H\lp f(p), f(q)\rp.
\end{equation} 
On the other hand, the Kobayashi pseudometric $K_M$ is the largest pseudometric so that all holomorphic maps $\H\ra M$
are nonexpanding. Thus, the Kobayashi pseudometric is bounded above by
$$
\delta_M(p,q) = \inf_{f\in \mathcal{O}(\H, M)} d_\H \lp f^{-1}(p), f^{-1}(q)\rp, 
$$
and if $\delta_M$ happens to satisfy the triangle inequality, then $K_M = \delta_M$.
In general,
\begin{equation}\label{D2}
K_M(p,q) = \inf \sum_{j=1}^{n} \delta_M(p_{j-1}, p_j),
\end{equation}
where the infimum is taken over all sequences $p_0,\ldots, p_n$ with $p_0 = p$ and $p_n = q$.
The Schwarz-Pick lemma implies that the assignments $M\mapsto K_M$ and $M\mapsto C_M$ satisfy condition
(i), while condition (ii) is a formal consequence of definitions \eqref{D1} and \eqref{D2}. In case $M$ is biholomorphic to a bounded domain in $\C^k$, the pseudometrics $K_M, C_M$ are nondegenerate and are thus referred to as the Carath\'eodory and Kobayashi metrics, respectively.

From the definitions, we see that every Schwarz-Pick system $d$ satisfies
$$
C_M \leq d_M \leq K_M
$$
for all complex manifolds $M$. Thus, if $C_M = K_M$, then every Schwarz-Pick system assigns
to $M$ the same pseudometric. The problem of determining for which manifolds $M$ the two pseudometrics
agree has attracted a great deal of attention.

\subsection{Complex geodesics and holomorphic retracts}
Let $d$ be a Schwarz-Pick system and $M$ a complex manifold. A {\em complex geodesic} for $d_M$ is a holomorphic and isometric embedding
$(\H, d_\H) \ra (M, d_M)$. A holomorphic map $\tau\in \mathcal{O}(\H, M)$ is said to be a {\em holomorphic retract} of $M$
if there exists a map $F\in \mathcal{O}(M,\H)$ so that $F\circ\tau = \id_\H$. We also say that $\tau$ {\em admits a holomorphic retraction}.
The main point of the following well-known lemma is that $\tau$ is complex geodesic for the Carath\'eodory metric if and only if it admits a holomorphic retraction \cite{Ja}.

\begin{Lemma}\label{retract}
Let $\tau:\H \ra M$ be a holomorphic map into a connected complex manifold. 
The following are equivalent:
\begin{enumerate}[label=\normalfont(\alph*), ref=(\alph*)]
\item
There is a pair of distinct $z,w \in \H$ so that $C_M\lp \tau(z), \tau(w) \rp = d_\H(z,w)$. \label{Sa}
\item
$\tau$ is a holomorphic retract of $M$. \label{Sb}
\item
$\tau$ is a complex geodesic for $C_M$. \label{Sc}
\item
$\tau$ is a complex geodesic for $K_M$ and the restrictions
of $C_M$ and $K_M$ to $\tau(\H)$ coincide. \label{Sd}
\end{enumerate}
\end{Lemma}

{\em Proof:}

\Implies{Sa}{Sb}:
There is a sequence of holomorphic maps $F_j : M \ra \H$
with $d_\H( F_j\circ \tau (z), F_j\circ \tau(w))$ converging to $d_\H(z, w)$. Postcomposing each $F_j$ by a M\"obius transformation,
we may assume $F_j\circ \tau$ fixes $z$ and maps $w$ to a point on the geodesic segment connecting $z$ and $w$. Then the $F_j$ form a normal
family. Any subsequential limit $F$ of the $F_j$ satisfies $F\circ \tau(z)=z$ and $F\circ \tau(w) = w$. By the Schwarz-Pick lemma, $F\circ \tau = \id_\H$.

\Implies{Sb}{Sc}:
Suppose $F:M \ra \H$ is holomorphic and satisfies $F\circ \tau = \id_\H$.
Then for any pair of points $z$ and $w$ in $\H$,$$d_\H(z,w) = d_\H( F\circ \tau(z), F\circ \tau(w)  ) \leq C_M( \tau(z), \tau(w) ) $$
because $F$ is holomorphic.
Also,$$C_M\lp \tau(z), \tau(w)\rp \leq d_\H(z,w)$$ because $\tau$ is holomorphic. 

\Implies{Sc}{Sd}:
For any $z,w$ in $\H$,
$$d_\H(z,w) = C_M\lp \tau(z),\tau(w)\rp \leq K_M(\tau(z),\tau(w)).$$
Since holomorphic maps decrease Kobayashi distance, the inequality must be an equality. 

\Implies{Sd}{Sa}:
Obvious.

\qed

{\em Remark:}
Let $p$ be a point in a connected complex manifold $M$.
To prove the implication \Implies{Sa}{Sb}, we used the fact that a family $\{F_j\}$ of holomorphic maps $M \ra \H$
is precompact in $\mathcal{O}(M,\H)$ if and only if $\{ F_j(p)\}$ is precompact in $\H$. (This is essentially a rephrasing of Montel's theorem.) We will use this fact throughout the paper.

\subsection{Symmetric spaces vs. Teichm\"uller space}
%\subsection{Hermitian symmetric spaces}
In case $M$ is a Hermitian symmetric space, the Kobayashi and Carath\'eodory metrics
coincide. Indeed, we have 
\begin{enumerate}
\item
Each pair of points in $M$ is contained in the image of a complex geodesic for $K_M$.
\item 
Every complex geodesic for $K_M$ is a holomorphic retract of $M$.
\end{enumerate}
So from Lemma \ref{retract}, we get $K_M = C_M$.
In fact, by a theorem of Lempert (\cite{Le1};\cite{Ja} Chapter 11) %\cite{Le1} \cite{Le2}
, the Kobayashi and Carath\'eodory metrics coincide for all bounded convex domains. 
(Every Hermitian symmetric space is biholomorphic to a bounded convex domain.)
%(See \cite{Ja} Chapter 11 for an expository treatment.) 
Whether the two metrics agree for all bounded $\C$-convex domains is an open question.

Given the many parallels between Teichm\"uller spaces and symmetric spaces, it is natural to ask whether the Carath\'eodory metric on Teichm\"uller space
is the same as Teichm\"uller-Kobayashi metric. As is the case with Hermitian symmetric spaces, any pair of points in $\T$ is contained in a complex geodesic for $K_\T$. The problem of determining whether the Kobayashi and Carath\'eodory metrics agree on
Teichm\"uller space thus reduces to checking whether each complex geodesic for $K_\T$ is a holomorphic retract.

%\begin{comment}
%\subsection{Symmetric spaces vs. Teichm\"uller spaces}
%In case $M$ is a Hermitian symmetric space of noncompact type, the Kobayashi and Carath\'eodory pseudometrics
%coincide. This follows from two facts:
%\begin{enumerate}
%\item
%Each pair of points $p,q \in M$ is contained in the image of a complex geodesic for $K_M$.
%\item 
%Every complex geodesic for $K_M$ is a holomorphic retract of $M$.
%\end{enumerate}

%In fact, by a theorem of Lempert (see \cite{Ja} Chapter 11) %\cite{Le1} \cite{Le2}
%, the metrics $K_M$ and $C_M$ coincide for all bounded convex domains.
%(See \cite{Ja} Chapter 11 for an expository treatment.) 
%Whether the two metrics agree for all bounded $\C$-convex domains is an open question.
%\end{comment}

%\subsection{The Kobayashi and Carath\'eodory Metrics on Teichm\"uller Space}
%Given the many parallels between Teichm\"uller spaces and symmetric spaces, it is natural to ask whether the Carath\'eodory metric on Teichm\"uller space
%is the same as Teichm\"uller-Kobayashi metric. As is the case with Hermitian symmetric spaces, the Teichm\"uller space $\T$ of a finite-type orientable surface embeds biholomorphically as a bounded domain
%in $\C^k$. Moreover, any pair of distinct points in $\T$ is contained in a unique complex geodesic $\tau$ for $K_\T$.
%Thus, the problem of determining whether the Kobayashi and Carath\'eodory metrics agree on
%Teichm\"uller space reduces to checking whether each complex geodesic for $K_\T$ is a retract.
 
\subsection{Abelian Teichm\"uller disks}
Complex geodesics for the Teichm\"uller-Kobayashi metric $K_\T$ are called {\em Teichm\"uller disks}. A Teichm\"uller disk 
is determined by the initial data of a point in $\T$ and a unit cotangent vector at that point. In other words, the disk is determined by
a unit-norm holomorphic quadratic differential $\phi$. We say $\phi$ {\em generates} the Teichm\"uller disk $\tau^{\phi}$.
(See Section \ref{T}.) If $\phi$ is the square of an Abelian differential, $\tau^{\phi}$ is called an {\em Abelian Teichm\"uller disk}.

Kra \cite{Kr} showed that the Kobayashi and Carath\'eodory metrics agree on certain subsets of Teichm\"uller space.
Namely, he proved
\begin{Theorem}\label{Kra} \cite{Kr}
Let $\T$ be the Teichm\"uller space of a finite-type orientable surface. If $\phi$ is a quadratic differential with no odd-order zeros,
then the restrictions of the metrics $K_\T$ and $C_\T$ to $\tau^\phi(\H)$ coincide. 
That is, $\tau^{\phi}$ is a complex geodesic for $C_\T$ and thus a holomorphic retract of $\T$.  
\end{Theorem}     
The key tool in the proof is the Torelli map from the Teichm\"uller space $\T_g$ of a closed surface to 
the Siegel upper half-space $\mathcal{Z}_g$. Kra showed that the Torelli map sends every Abelian Teichm\"uller disk in $\T_g$
to a complex geodesic in the symmetric space $\mathcal{Z}_g$. Post-composing by a holomorphic retraction $\mathcal{Z}_g \ra \H$
onto this complex geodesic yields a retraction $\T_g \ra \H$ onto the Abelian Teichm\"uller disk.
A covering argument then extends the result to all differentials $\phi$ with no odd-order zeros. (Note that, at a puncture, $\phi$ can have a simple pole
or a zero of any order, and the theorem still holds.)

\subsection{Carath\'eodory $\neq$ Teichm\"uller}
However, in \cite{Mar} it was shown that the Carath\'eodory and Kobayashi metrics on Teichm\"uller space do not coincide:
\begin{Theorem}\label{neq} \cite{Mar}
The Kobayashi and Carath\'eodory metrics on the Teichm\"uller space of a closed surface of genus at least two do not coincide;
if $g\geq 2$, there is a Teichm\"uller disk in $\T_g$ which is not a holomorphic retract.  
\end{Theorem}

An elementary covering argument, outlined in the appendix of this paper, extends the result
to all Teichm\"uller spaces $\T_{g,n}$ of complex dimension at least two:
\begin{Theorem}\label{neq2}
Suppose $\dim_\C \T_{g,n} := 3g-3+n \geq 2$.
The Kobayashi and Carath\'eodory metrics on $\T_{g,n}$ are different.
\end{Theorem}

{\em Remark:} 
The one-dimensional spaces $\T_{0,4}$ and $\T_{1,1}$ are biholomorphic
to $\H$, so the Kobayashi and Carath\'eodory metrics are equal to the Poincar\'e metric.

Theorem \ref{neq2} has consequences for the global geometry of Teichm\"uller space.
Teichm\"uller space is biholomorphic via Bers' embedding to a bounded domain in $\C^k$.
However, combined with Lempert's theorem, Theorem \ref{neq2} implies
\begin{Theorem}
The Teichm\"uller space $\T_{g,n}$ is not biholomorphic to a bounded convex domain in $\C^k$,
whenever $\dim_\C \T_{g,n}\geq 2$.
\end{Theorem}
See \cite{Gu} for related convexity results. See \cite{An1} and \cite{An2} for other recent results comparing the complex
geometry of Teichm\"uller spaces and symmetric spaces.

\subsection{Outline}
It remains to characterize the quadratic differentials which generate holomorphic retracts. 
Our Conjecture \ref{Conj} is that the converse of Kra's result holds $-$ $\tau^{\phi}$ is a holomorphic retract if and only if $\phi$ has no odd-order zeros.
%In other words, the conjecture is that, if $\phi$ has an odd-order zero, then $\tau^{\phi}$ is not a holomorphic retract.

In the rest of the paper, we develop some tools towards a proof of Conjecture \ref{Conj}. The conjecture is obviously true for the Teichm\"uller
spaces $\T_{1,1}$ and $\T_{0,4}$ of complex dimension one.
Our {\bf main result}, Theorem \ref{Main}, is that the conjecture holds for the spaces $\T_{0,5}$ and $\T_{1,2}$ of complex dimension two.

The idea of the proof is as follows. 
We first show that the property of generating a holomorphic retract is a closed condition on the bundle $\widetilde{\mathcal{Q}}$ of marked quadratic differentials.
The condition is also invariant under the actions of $\SL_2(\R)$ and the mapping class group.
Thus, it descends to a closed, $\SL_2(\R)$-invariant
condition on the moduli space $\mathcal{Q}$ of unmarked quadratic differentials. 
In other words, if $\phi$ generates holomorphic retract, then so does every element
of its $\SL_2(\R)$ orbit closure in $\mathcal{Q}$. 
On the other hand, the $\SL_2(\R)$-orbit closure of any quadratic differential contains a Jenkins-Strebel differential in the same stratum \cite{Mi}\cite{Sm}.

To prove Conjecture \ref{Conj}, it thus suffices to establish that no Jenkins-Strebel differential with an odd-order zero generates a retract.
To this end, we prove an analytic criterion characterizing Jenkins-Strebel differentials which generate retracts.
Given a Jenkins-Strebel differential $\phi$ with $k$ cylinders, we define a holomorphic map $\mathcal{E}^{\phi}:\H^k \ra \T$
called a {\em Teichm\"uller polyplane}. The marked surface $\mathcal{E}^{\phi}(\lambda_1,\ldots, \lambda_k)$
is obtained by applying the map $x+iy\mapsto x+\lambda_jy$ to the $j$th cylinder (Figure \ref{Genus2Shear}). The Teichm\"uller disk $\tau^\phi$ is the
diagonal of $\mathcal{E}^{\phi}$, so if $F: \T \ra \H$ is a holomorphic retraction onto $\tau^{\phi}(\H)$, then
the composition $f=F\circ \tau : \H^k \ra \H$ restricts to the identity on the diagonal:
 $$f(\lambda,\ldots, \lambda) = \lambda.$$
Our analytic criterion, Theorem \ref{criterion}, states that if a unit area Jenkins-Strebel differential $\phi$ generates a retract,
then the retraction $F$ can be chosen so that $f = F\circ \mathcal{E}^{\phi}$ is a convex combination of the coordinate functions: 
$$f(\lambda_1,\ldots, \lambda_k) = \sum_j a_j \lambda_j,$$
where $a_j$ is the area of the $j$th cylinder. In other words, $\tau^{\phi}$ is a retract if and only if the linear
function $\sum_j a_j \lambda_j$ on the polyplane $\mathcal{E}^\phi(\H^k)$ extends to a holomorphic map
$\T \ra \H$. As a corollary of this criterion, we observe that if a Jenkins-Strebel differential $\phi$
generates a retract, then so does any differential obtained by horizontally shearing the cylinders of $\phi$.

To prove Conjecture \ref{Conj} for $\T_{0,5}$, let $\phi\in {\mathcal{Q}}_{0,5}$ be a Jenkins-Strebel
differential with an odd-order zero. So $\phi$ has a simple zero and five poles.
Using the results of \cite{Mi},\cite{Sm} and a simple combinatorial argument, we show that $\phi$ contains in its $\SL_2(\R)$ orbit closure a {\em two-cylinder} Jenkins-Strebel differential
$\phi'$ with a simple zero. Shearing the cylinders of $\phi'$ yields an {\em L-shaped pillowcase differential} $\psi$ (Figure \ref{Lshaped}).
Now, assume for the sake of contradiction that $\phi$ generates a retract. Then so does the differential $\phi'$ in $\ol{SL_2(\R)\phi}$,
and so does the L-shape $\psi$ obtained by shearing the cylinders of $\phi'$. But in \cite{Mar}, Markovic shows that an L-shape does not satisfy our criterion; there is no holomorphic map $F:\T_{0,5} \ra \H$ extending $a_1\lambda_1 + a_2\lambda_2$. The idea of Markovic's proof is to assume a holomorphic extension $F$ exists and then, using the Schwarz-Christoffel formula, obtain a contradiction on the smoothness of $F$ at the boundary of the bidisk $\mathcal{E}^\psi(\H^2)$.

This proves the Conjecture \ref{Conj} for the sphere with five punctures. The isomorphism $\T_{0,5} \cong \T_{1,2}$ yields the corresponding result for the twice-punctured torus.

\section{Dynamics on moduli space}\label{DYN}
In this section, we describe the role of dynamics in the classification of complex geodesics for the Carath\'eodory metric.
After recalling some basic definitions, we show how the $\GL_2^+(\R)$ action on $\mathcal{Q}$ allows
us to reduce Conjecture \ref{Conj} to the case of Jenkins-Strebel differentials. 
Using ergodicity of the Teichm\"uller geodesic flow, we show that most quadratic differentials do not generate holomorphic retracts.

\subsection{The $\GL_2^+$ action on $\widetilde{\mathcal{Q}}$}
\begin{figure}[h]
\includegraphics{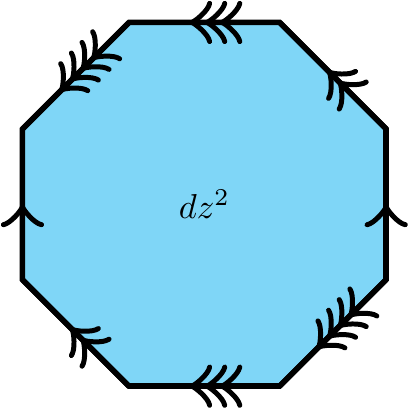}
\caption{A quadratic differential on the surface of genus two. The vertices glue up to a single cone point of angle $6\pi$,
corresponding to an order four zero of the differential.}
\label{Genus2}
\end{figure}

\begin{figure}[h]
\includegraphics{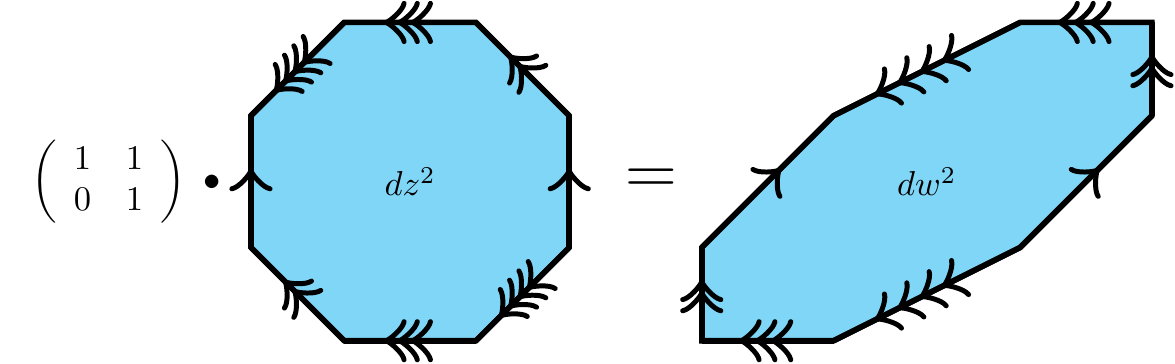}
\caption{The action of a shear on the differential in Figure \ref{Genus2}. The underlying surface of the resulting quadratic differential is 
$\tau^{\phi}(1+i)$.}
\label{Genus2Shear}
\end{figure}

Let $\T=\T_{g,n}$ be the Teichm\"uller space of marked complex structures on a finite-type, orientable surface $S_{g,n}$ of genus $g$ with $n$ punctures.
Let $\widetilde{\mathcal{Q}}=\widetilde{\mathcal{Q}}_{g,n}$ denote the bundle of marked, nonzero, integrable, holomorphic quadratic differentials over
$\T$. Equivalently, $\widetilde{\mathcal{Q}}$ is the bundle of marked half-translation structures on $S_{g,n}$ (Figure \ref{Genus2}).
The group $\GL_2^+$ of orientation-preserving
linear maps $\R^2 \ra \R^2$ acts on $\widetilde{\mathcal{Q}}$ by post-composition of flat charts.  In other words, the action is by affine deformations of the polygonal
decomposition of a differential (Figure \ref{Genus2Shear}).

\pagebreak

\subsection{Teichm\"uller disks}\label{T}
Let $p:\widetilde{\mathcal{Q}} \ra \T$ denote the projection sending a quadratic differential to its underlying Riemann surface. 
The action of a conformal linear transformation does not change the underlying Riemann surface of a quadratic differential.
Therefore, $$\tau^{\phi}:g\mapsto p(g\cdot \phi)$$ is a well-defined map from the upper half-plane $\H\cong \C^\times \setminus\GL^+_2$ to
Teichm\"uller space $\T$.
To give a more explicit description of the map $\tau^{\phi}$, note that $\H$ sits in $\GL^+_2$ as the subgroup generated by vertical stretches and horizontal shears:
$$ \left\{\lp\begin{array}{cc} 1 & \Re(\lambda) \\ 0 &\Im(\lambda) \end{array}\rp \middle| \lambda\in \H\right\}.$$
The Teichm\"uller disk generated by $\phi$ is
\begin{equation}\label{TWO}\lambda \mapsto p(\lambda\cdot \phi).
\end{equation}
Written in complex coordinates, the action of the matrix $\lp\begin{array}{cc} 1 & \Re(\lambda) \\ 0 &\Im(\lambda) \end{array}\rp$
is 
\begin{equation*}
x+iy \mapsto x+\lambda y,
\end{equation*}  
which has Beltrami coefficient
$$\frac{i-\lambda}{i+\lambda}\frac{\ol{dz}}{dz}.$$
Thus, $\tau^\phi(\lambda)$ is the quasiconformal deformation of $X:=p(\phi)$ with Beltrami coefficient
$$\frac{i-\lambda}{i+\lambda}\mu_\phi,$$
with $\mu_\phi := \phi^{-1}\lpi \phi \rpi$. 
The Teichm\"uller disk $\tau^\phi$ is the unique Kobayashi geodesic with initial data $\tau(i) = X$ and $\tau'(i) = \frac{i}{2}\mu_\phi$.

\subsection{Orbit closures}
The mapping class group $\text{MCG}$ of $S_{g,n}$ acts on $\T$ and $\widetilde{\mathcal{Q}}$ by changes of marking.
The action of each mapping class is a biholomorphism, and by a theorem of Royden \cite{Ro}, every biholomorphism
of Teichm\"uller space arises in this way (as long as $\dim_\C \T\geq 2$).
The quotient of Teichm\"uller space by the $\text{MCG}$ action is the moduli space
of complex structures on $S_{g,n}$. The quotient of $\widetilde{\mathcal{Q}}$ by the $\text{MCG}$ action
is the space $\mathcal{Q}$ of (unmarked) half-translation structures on $S_{g,n}$.
The $\GL_2^+$ action on $\widetilde{\mathcal{Q}}$ descends to an action on $\mathcal{Q}$.  

Let $\phi\in \widetilde{\mathcal{Q}}$ and $\alpha\in \text{MCG}$. Then the disk $\tau^{\phi}$ is a holomorphic retract if and only if
$\tau^{\alpha(\phi)}$ is. Indeed, if $F:\T \ra \H$ is a retraction onto $\tau^{\phi}$, then $F\circ \alpha^{-1}$ is a retraction
onto $\tau^{\alpha(\phi)} = \alpha\circ \tau^{\phi}$. Thus, we will say $\phi\in \mathcal{Q}$ generates a retract
if every element of its preimage in $\widetilde{\mathcal{Q}}$ does.
The property of generating a retract is also invariant under the $\GL_2^+$ action. Indeed, if 
$$g=\lp \begin{array}{cc} a & b \\ c & d \end{array} \rp\in \GL_2^+,$$ and $m:\H\ra \H$ is the associated M\"obius
transformation
$$m(\lambda) = \frac{d\lambda+b}{c\lambda+a},$$
then 
$\tau^{g\cdot \phi} = \tau^{\phi} \circ m$. 

The following Proposition is key:
\begin{Proposition}\label{Closure}
If $\phi$ generates a holomorphic retract, so does every element in the orbit closure $\ol{\GL_2^+\phi}\subset \mathcal{Q}$. 
\end{Proposition}
{\em Proof:}
From the above discussion, we know that every element of the orbit $\GL_2^+\phi$ generates a retract.
The desired result will follow from the next Lemma.

\begin{Lemma}
Suppose the sequence $\phi^1,\phi^2,\ldots$ converges to $\phi$ in $\widetilde{\mathcal{Q}}$.
Let $\tau^N$, $\tau$ be the Teichm\"uller disks generated by $\phi^N$ and $\phi$, respectively.
If each $\tau^{N}$ is a holomorphic retract, then so is $\tau$.
\end{Lemma}
{\em Proof:} This follows by the continuity of the Kobayashi and Carath\'eodory metrics on $\T$ (see \cite{Ea}).
However, we prefer to give a direct proof.

For each $N$, let $F^N: \T \ra \H$ be a holomorphic map satisfying $F^N\circ \tau^{N} = \id_\H$. 
Then $\{F^N\}$ is a normal family. To see this, let $X^N$ and $X$ denote the marked surfaces $p(\phi^N), p(\phi)$.
Then
$$d_\H\lp  F^N(X), i    \rp = d_\H\lp  F^N(X), F^N(X^N)    \rp \leq K_\H(X, X^N),$$
which is uniformly bounded in $N$.
Thus, passing to a subsequence, we may assume $F^N$ converges locally uniformly to a holomorphic map $F:\T \ra \H$.
By continuity of the $\GL_2^+$ action on $\widetilde{\mathcal{Q}}$, the sequence $\tau^{N}$ converges locally uniformly to $\tau$.
Therefore, 
$$\Phi\circ \tau^{\phi} = \lim_{N\ra \infty} \Phi^N\circ \tau^{N} = \id_\H.$$ 
\qed

Let $\mathcal{Q}_1$ denote the space of unit-area half-translation surfaces.
The {\em Teichm\"uller geodesic flow} on $\mathcal{Q}_1$  is the action of the subgroup
$$
\left\{ \lp \begin{array}{cc} e^t & 0\\ 0 & e^{-t} \end{array}\rp \middle |~ t\in \R \right\} \subset \GL_2^+. 
$$ 
With respect to a suitable probability measure on $\mathcal{Q}_1$, the Teichm\"uller geodesic flow is ergodic \cite{Mas}, \cite{Ve}.
In particular, for almost all $\phi \in \mathcal{Q}_1$, the orbit $\GL_2^+ \phi$ is dense in $\mathcal{Q}_1$.
Combined with Proposition \ref{Closure} and Theorem \ref{neq2}, this implies
\begin{Theorem}\label{Most}
For almost every quadratic differential $\phi \in \mathcal{Q}_1$, the Teichm\"uller disk
$\tau^{\phi}(\H)$ is not a holomorphic retract.
\end{Theorem}

\subsection{The horocycle flow and Jenkins-Strebel differentials}\label{Horocycle}
Recall that $\phi\in \mathcal{Q}$ is said to be {\em Jenkins-Strebel} if its nonsingular horizontal
trajectories are compact.

The {\em horocycle flow} on $\mathcal{Q}$ is the action of the subgroup
$$
H=\left\{ \lp \begin{array}{cc} 1 & t\\ 0 & 1 \end{array}\rp \middle |~ t\in \R \right\} \subset \GL_2^+. 
$$
Minsky and Weiss \cite{Mi} showed that every closed $H$-invariant set of $\mathcal{Q}_1$
contains a minimal closed $H$-invariant subset. If $\mathcal{Q}_1$ were compact, this would 
follow by a standard Zorn's lemma argument. The main point of their proof is a {\em quantitative nondivergence} result
which states that each $H$-orbit spends a large fraction of its time in the $\vn$-thick part of $\mathcal{Q}_1$.  

Subsequently, Smillie and Weiss \cite{Sm}
showed that every minimal closed $H$-invariant set is the $H$ orbit closure of a Jenkins-Strebel differential.
In particular, the orbit closure $\ol{H\phi}$ of any $\phi \in \mathcal{Q}$ contains a Jenkins-Strebel differential.
As Smillie and Weiss observed, the above results continue to hold  for the action of $H$ on each stratum of $\mathcal{Q}$.
In particular, this means that if $\phi$ has an odd-order zero, then $\ol{H\phi}$
contains a Jenkins-Strebel differential with an odd-order zero. 

An advantage of working with the horocycle flow, rather than the full $\GL_2^+$ action,
is that $H$ preserves horizontal cylinders. That is, if $\phi$ has a horizontal cylinder, then every element of $\ol{H\phi}$
has a cylinder of the same height and length. Suppose in addition that the cylinder of $\phi$ is not dense in $S_{g,n}$.
Then since $\| h\cdot\phi\| = \| \phi\|$ for all $h$ in $H$, a Jenkins-Strebel differential in $\ol{H\phi}$ has at least two cylinders.

We summarize the above in the following theorem.
\begin{Theorem}\label{hclosure}
Let $\phi \in \mathcal{Q}$ be a quadratic differential.
\begin{enumerate}[label=\normalfont(\alph*), ref=(\alph*)]
\item
The closure $\ol{H\cdot\phi}$ contains a Jenkins-Strebel differential $\psi$.
If $\phi$ has an odd-order zero, then $\psi$ can be taken to also have an odd-order zero.
\item
If $\phi$ has a horizontal cylinder which is not dense in $S_{g,n}$ then $\psi$
has at least two cylinders.
\end{enumerate} 
\end{Theorem}

By Proposition \ref{Closure} and Theorem \ref{hclosure}(a), our Conjecture \ref{Conj} reduces to the case of Jenkins-Strebel differentials.
That is, it suffices to show that no Jenkins-Strebel differential with an odd-order zero generates a holomorphic retract.

%\begin{figure}[h]
%\footnotesize
%\stackunder[5pt]{\includegraphics[scale=.65]{Singularity-1}}{Pole at a singularity}
%\hspace{1cm}%
%\stackunder[5pt]{\includegraphics[scale=.65]{Singularity0}}{Regular point}
%\hspace{1cm}%
%\stackunder[5pt]{\includegraphics[scale=.65]{Singularity1}}{Order 1 zero}
%\hspace{1cm}%
%\stackunder[5pt]{\includegraphics[scale=.65]{Singularity2}}{Order 2 zero}
%\caption{Local structure of the horizontal foliation}
%\label{Singularities}
%\end{figure}

\section{Jenkins-Strebel differentials and Teichm\"uller polyplanes}\label{POLY}
In this section, we begin our analysis of Teichm\"uller disks generated by Jenkins-Strebel differentials.
The key observation is that such a disk is the diagonal of a certain naturally defined polydisk holomorphically embedded in $\T$. 

\subsection{Teichm\"uller polyplanes}
\begin{figure}[h]
\includegraphics{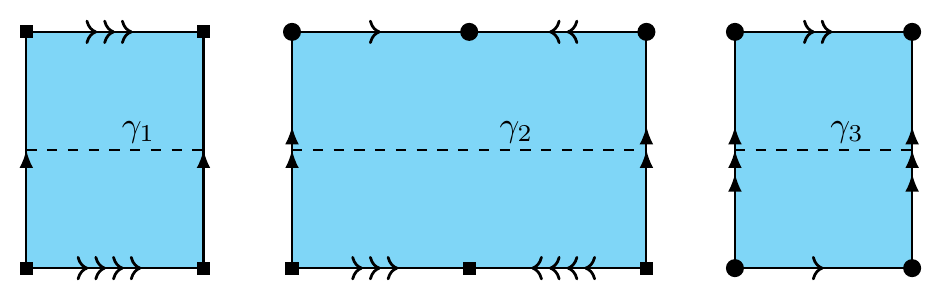}
\caption{A Jenkins-Strebel differential on a genus 2 surface. The differential has two order 2 zeros, indicated by the dot and the square.}
\label{StrebelGenus2}
\end{figure}

\begin{figure}[h]
\includegraphics{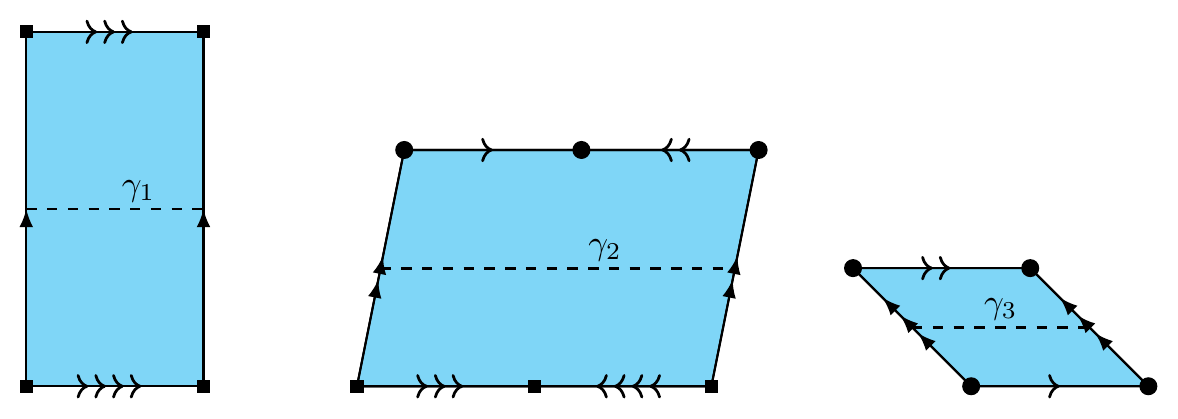}
\caption{The action of an element of $\H^3$ on the differential from Figure \ref{StrebelGenus2}. The resulting Riemann surface
is $\E^{\phi}(1.5i,.2+i,  -.5 + .5i)$}
\label{StrebelGenus2Shear}
\end{figure}

\begin{figure}[h]
\includegraphics{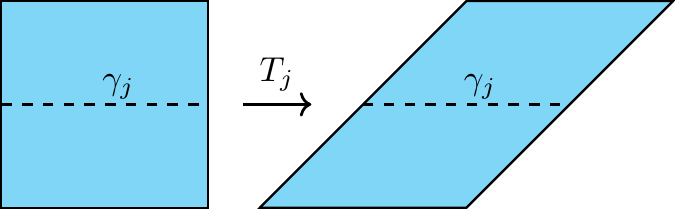}
\caption{Translation by $m_j^{-1}$ in $\H^k$ corresponds to a Dehn twist $T_j$ about $\gamma_j$ in the Teichm\"uller space.}
\label{DehnTwist}
\end{figure}

The core curves of a Jenkins-Strebel differential form a collection of essential simple closed curves, which are pairwise disjoint and non-homotopic
(Figure \ref{StrebelGenus2}).
We will call such a collection of curves a {\em disjoint curve system}.
Let $\mathcal{C} = \{ \gamma_1,\ldots , \gamma_k\}$ be a disjoint curve system on $S_{g,n}$.
We say a Jenkins-Strebel differential is of type $\mathcal{C}$ if the cores of its cylinders are homotopic to $\gamma_1,\ldots, \gamma_k$.
We define an action of the $k$-fold product $\H^k=\H\times \cdots \times \H$ on the differentials of type $\mathcal{C}$.
The tuple $\llambda = (\lambda_1,\ldots, \lambda_k)\in \H^k$ acts on the $j$th cylinder of $\phi$ by the affine map $x+iy \mapsto x+\lambda_j y$ 
(Figure \ref{StrebelGenus2Shear}). Since this map takes horizontal circles isometrically to horizontal circles, the result is a well-defined Jenkins-Strebel
differential $\llambda\cdot \phi$.
Projecting the orbit of $\phi$ to the Teichm\"uller space, we get a map $\E^{\phi}: \H^k \ra \T$ defined by
$$\E^{\phi}(\llambda) = p(\llambda\cdot \phi).$$
We call $\E^{\phi}$ the {\em Teichm\"uller polyplane} associated to $\phi$. Below, we list some properties of Teichm\"uller polyplanes.
\begin{enumerate}[label={}, leftmargin=0cm]
\item
{\bf The Teichm\"uller disk associated to $\phi$ is the diagonal of the polyplane:}
$$\tau^{\phi}(\lambda) = \E^{\phi}(\lambda,\ldots, \lambda).$$
\item
{\bf The Teichm\"uller polyplane mapping sends translations to Dehn twists} (See Figure \ref{DehnTwist}.)
To make this precise, let $m_j$ denote the modulus (height divided by length) of
the $j$th cylinder. Let $T_j$ denote the Dehn twist about the core curve $\gamma_j$. Then
\begin{equation}\label{Twist}
\E^{\phi}\lp \llambda  + (0,\ldots, m_j^{-1}, \ldots, 0) \rp = T_j\circ \E^{\phi}(\llambda).
\end{equation}
Equation \eqref{Twist} is crucial. It will allow us to relate the analysis of holomorphic maps $\H^k \ra \H$
to the geometry of Teichm\"uller space.
\item
{\bf The polyplane mapping is a holomorphic embedding.}
The mapping is holomorphic because the Beltrami coefficient of $x+iy\mapsto x+\lambda_j y$ is holomorphic in $\lambda_j$.
Holomorphicity implies that $\E^{\phi}$ is nonexpanding for the Kobayashi metrics on $\H^k$ and $\T$:
\begin{equation}\label{Schwarz}
K_\T\lp \E^{\phi}(\llambda_1), \E^{\phi}(\llambda_2)\rp \leq K_{\H^k}(\llambda_1,\llambda_2)
\end{equation}
(Recall that the Kobayashi metric $K_{\H^k}$ is the supremum of the Poincar\'e metrics on the factors.)
We prove that $\E^{\phi}$ is an embedding in Theorem \ref{Embed} below.  
\item
{\bf However, the mapping is not proper.}
Indeed, there are no proper holomorphic maps $\H^k \ra \T$. (See Page 75, Corollary 1 of \cite{Ta}.)
As one of the one of the coordinate functions approaches the real axis and the height of the corresponding cylinder goes to 0,
the sequence of image points may converge in $\T$. This lack of properness was critical in the proof \cite{Mar} that the Kobayashi and Carath\'eodory
metrics on $\T$ are different.
\end{enumerate}
  
\pagebreak
  
\subsection{Teichm\"uller polyplanes are embedded}\label{EMBED}
Below, we prove that the Teichm\"uller polyplane $\E^{\phi}:\H^k \ra \T$ is an embedding.
We do not use this fact in the proof of our main result; the reader may choose to skip this section. 

Let $\phi$ be a quadratic differential and $\gamma$ a closed curve. We denote by $L_\phi(\gamma)$ the $\phi$-length of $\gamma$ $-$ the shortest length of a curve homotopic to $\gamma$,
measured in the flat metric associated to $\phi$. If $\gamma$ is the core of a cylinder of $\phi$, then $L_{\phi}(\gamma)$
is the circumference of the cylinder.
The following result of Jenkins \cite{Je} asserts uniqueness of
Jenkins-Strebel differential with given length data.
\begin{Proposition}\label{Jenk}
Let $\mathcal{C}=\{\gamma_1,\ldots, \gamma_k\}$ be a disjoint curve system,
and let $l_1, \ldots, l_k$ be positive numbers. There is at most one Jenkins-Strebel differential $\phi$ whose
core curves are homotopic to a subset of $\mathcal{C}$ and which satisfies $L_{\phi}(\gamma_j) = l_j$
for $j=1,\ldots, k$.
\end{Proposition}
Combined with Proposition \ref{Jenk}, the following result implies injectivity of $\E^{\phi}$.
\begin{Lemma}\label{free}
The action of $\H^k$ on Jenkins-Strebel differentials of given type $\{ \gamma_1,\ldots, \gamma_k \}$ is free.
In other words, each orbit map $\llambda\mapsto \llambda\cdot \phi$ is injective.
\end{Lemma}
{\em Proof:} Suppose $\llambda \cdot \phi = \mmu \cdot \phi$. Then the height of the $j$th cylinder is the same
for $\llambda \cdot \phi$ and $\mmu \cdot \phi$, so $\Im (\lambda_j) = \Im (\mu_j)$.  

We need to show $\Re(\lambda_j) = \Re(\mu_j)$. Suppose not.
Then $\mmu = \llambda + \t$ for some nonzero vector $\t \in \R^k$. The equation 
$\lp \llambda + \t\rp \cdot \mu = \llambda \cdot \mu$ combined with the fact that $(\llambda, \phi)\mapsto \llambda\cdot \phi$
is a group action implies 
$$\lp \llambda + N\t\rp \cdot \phi = \llambda \cdot \phi.$$
for every positive integer $N$.
Projecting to Teichm\"uller space, we get
$$\mathcal{E}^{\phi}(\llambda+N\t) = \mathcal{E}^{\phi}(\llambda).$$
But this is impossible since for large $N$, $\mathcal{E}^{\phi}(\llambda+N\t)$ 
is bounded distance from a translate of $\mathcal{E}^{\phi}(\llambda)$ by a big Dehn multi-twist.

To make the argument precise, let $\v^N$ denote the vector with $j$th entry 
$$\frac{\lfloor Nt_jm_j  \rfloor}{ m_j}.$$
(The vector $\v^N$ is an element of the lattice $\bigoplus m_j^{-1}\Z$ which approximates $N\t$.)
Now by \eqref{Twist}, $\mathcal{E}^{\phi}(\llambda + \v^N)$ is the marked surface obtained by twisting
$\lfloor Nt_jm_j  \rfloor$ times about $\gamma_j$. By proper discontinuity of the action of the mapping class group on $\T$,
the sequence $\mathcal{E}^{\phi}(\llambda + \v^N)$ leaves every compact set as $N\ra\infty$.
However, by \eqref{Schwarz},
\begin{align*}
K_{\T}\lp \mathcal{E}^{\phi}(\llambda + \v^N), \mathcal{E}^{\phi}(\llambda) \rp &= 
K_{\T}\lp \mathcal{E}^{\phi}(\llambda + \v^N), \mathcal{E}^{\phi}(\llambda + N\t) \rp \\
&\leq K_{\H^k}(\llambda + \v^N, \llambda + N\t),
\end{align*}
which is bounded by a constant independent of $N$. This is a contradiction.
\qed

We now come to the main result of this section.
\begin{Theorem}\label{Embed}
The Teichm\"uller polyplane
$\mathcal{E}^{\phi}:\H^k \ra \T$ is a holomorphic embedding.
\end{Theorem}
{\em Proof:} We have already seen that $\mathcal{E}^{\phi}$ is holomorphic.
To prove injectivity, suppose $\mathcal{E}^{\phi}(\llambda) = \mathcal{E}^{\phi}(\mmu)$.
Then $\llambda\cdot \phi$ and $\mmu\cdot \phi$ are quadratic differentials on the same marked Riemann surface.
By construction, the cylinders of these two differentials have the same lengths. So  $\llambda\cdot \phi = \mmu\cdot \phi$
by Proposition \ref{Jenk}. Now, by Lemma \ref{free}, $\llambda = \mmu$.

It remains to show that $\E^{\phi}$ is a homeomorphism onto its image.
To this end, let $\llambda^1, \llambda^2, \ldots$ be a sequence which leaves all compact subsets of $\H^k$.
We must verify that $\E^{\phi}(\llambda^1), \E^{\phi}(\llambda^2), \ldots$ does not converge to an element
of $\E(\H^k)$. It suffices to check the following cases:
\begin{enumerate}[label=\bf(\roman*), leftmargin=*]
\item
{\bf The imaginary part of some component of $\llambda^N$ goes to infinity as $N\ra\infty$.}

\noindent Then the modulus of the corresponding cylinder goes to infinity. Therefore,
the extremal length (see e.g. \cite{Ke}) of the core curve converges to 0. Hence, 
$\E^{\phi}(\llambda^N)$ leaves all compact subsets of $\T$.

\item
{\bf The imaginary parts of all components $\llambda^N$ are bounded above but at least one of them converges to zero.}

\noindent Suppose $\E^{\phi}(\llambda^N)$ converges to a marked surface $X$.
Since the imaginary parts of $\llambda^N$ stay bounded, the norms
$\| \llambda^N\cdot \phi\|$ stay bounded. Passing to a subsequence,
we may thus assume that $\llambda^N\cdot \phi$
converges to a Jenkins-Strebel differential $\psi$ on $X$. By continuity, $\psi$ has the same length 
data as every element of the orbit $\H^k\cdot \phi$, that is $L_\psi(\gamma_j) = L_\phi(\gamma_j)$ for $j=1,\ldots, k$.
However, if the $j$th coordinate $\lambda^N_j$ converges to 0, then the $j$th cylinder of $\psi$ is degenerate.
In other words, $\gamma_j$ is not a core curve of $\psi$. So $\psi$ is not in $\H^k\cdot \phi$. By Proposition \ref{Jenk}, $X$ is not in $\mathcal{E}^\phi(\H^k)$.

\item
{\bf All of the imaginary parts stay bounded between two positive numbers,
but some of the real parts go to infinity.}

\noindent Arguing as in the proof of Lemma \ref{free}, there is a sequence of multi-twists $\alpha^N:\T\ra\T$ so that $$K_\T\lp \E({\llambda^N}) , \alpha^N \circ \E^{\phi}(i) \rp$$
is bounded uniformly in $N$. Since $\alpha^N \circ \E^{\phi}(i)$ leaves all compact subsets of $\T$,
so does $\E({\llambda^N})$.  \qed   
\end{enumerate}
Case (ii) accounts for the fact that $\mathcal{E}^{\phi}$ is not proper.

\section{The analytic criterion}\label{CRIT}
The goal of this section is to prove the following analytic criterion characterizing Jenkin-Strebel differentials which generate retracts. 
The criterion generalizes results from \cite{Ge} \cite{Mar}.

\begin{Theorem}\label{criterion}
Let $\mathcal{\phi}\in \widetilde{\mathcal{Q}}$ be a unit-area Jenkins-Strebel differential, and let $a_j$ denote the area of the $j$th cylinder of $\phi$.
Let $\mathcal{E}^\phi:\H^k \ra \T$ be the Teichm\"uller polyplane associated to $\phi$.
Then the Teichm\"uller disk $\tau^{\phi}$ is a holomorphic retract if and only if there exists a holomorphic map $G: \T \ra \H$ so that
\begin{equation}\label{DUH}
G\circ \E^\phi (\llambda) = \sum_{j=1}^k a_j \lambda_j.  
\end{equation}
\end{Theorem}

In other words, $\tau^{\phi}$ is a retract if and only if the function $\sum_{j=1}^k a_j \lambda_j$ on $\mathcal{E}^{\phi}(\H^k)$ admits a
holomorphic extension to the entire Teichm\"uller space. Heuristically, the more cylinders $\phi$ has, the stronger the criterion is.
If $\phi$ has one cylinder, the criterion is vacuous. If the core curves of $\phi$ form a maximal disjoint curve system,
then the polyplane is an open submanifold of $\T$, so the criterion says that $\sum_{j=1}^k a_j \lambda_j$ has a unique extension to
a holomorphic map $\T \ra \H$.

{\em Remark 1:} Alex Wright has pointed out to us a proof, based on his work in \cite{Wr}, that the orbit closure of a differential $\phi$ with an odd-order
zero contains a Strebel differential with at least two cylinders,
except potentially if $\phi$ is a pillowcase cover. So the criterion in Theorem \ref{criterion} almost always gives us at least some nontrivial information.

{\em Remark 2:} A potential program to prove Conjecture \ref{Conj} is to
\begin{enumerate}[label=(\roman*)]
\item
Use the criterion to identify a class of Jenkins-Strebel differentials which do not generate retracts.
\item
Show that the orbit closure of any differential with an odd-order zero contains a Jenkins-Strebel differential of that class.
\end{enumerate}
We will carry out this program for $\T_{0,5}$ by working with the class of differentials with two cylinders.

We return to the proof of Theorem \ref{criterion}. 
The ``if" direction is easy; if there is a holomorphic map $G:\T \ra \H$
satisfying \eqref{DUH}, then
\begin{align*}
G\circ \tau^{\phi}(\lambda) &= 
G\circ \mathcal{E}^{\phi} (\lambda,\ldots, \lambda)\\ 
&=  \sum_{j=1}^k a_j \lambda = \lambda.
\end{align*}
Thus, $G$ is a holomorphic retraction onto $\tau^{\phi}(\H)$.

 To prove the other direction, suppose there is a map
$F:\H \ra \T$ so that $F\circ \tau^{\phi} = \id_\H$. The idea of the proof is to approximate the desired $G$ by maps
of the form $t + F\circ \alpha$, with $t\in \R$ a translation and $\alpha \in \MCG$ a multi-twist. First, we recall a lemma from \cite{Mar}.

\begin{Lemma}
Let $f:\H^k \ra \H$ be the composition $F\circ\E^\phi$. Then $f$ satisfies
\begin{equation}\label{Diag}
f(\lambda,\ldots, \lambda)=\lambda
\end{equation}
and
\begin{equation}\label{Deriv}
\pd{f}{\lambda_j}(i,\ldots, i)=a_j,
\end{equation} 
where $a_j$ is the area of the $j$th cylinder
of the unit-area differential $\phi$. 
\end{Lemma}
{\em Proof: }
Equation $\eqref{Diag}$ is a restatement of $F\circ \tau^{\phi}=\id_\H$.
To prove \eqref{Deriv}, let $X=\tau^{\phi}(i)$ be the underlying surface of $\phi$. We first show that the cotangent vector $dF_X \in T^*_X\T$ is represented by the quadratic differential $-2i\phi$.
To this end, let $\mu_{\phi}$ be the Beltrami differential $\phi^{-1}\lpi\phi\rpi$. By definition, $\lp \tau^{\phi}\rp'(i) = \frac{i}{2}\mu_\phi$. So by the chain rule,
$$dF_X \lp \frac{i}{2}\mu_{\phi} \rp = 1.$$ But also $$\int_X (-2i\phi) \lp \frac{i}{2}\mu_\phi\rp=\|\phi\|=1.$$ 
Since $F$ is holomorphic, $dF_X$ has norm at most 1 with respect to the infinitesimal Kobayashi metrics on $T_X\T$ and $T_i\H$.
Since the infinitesimal Kobayashi metric on $T_i\H$ is half the Euclidean metric, the unit norm ball of $\Hom_\C(T_X\T, T_i\H)$
corresponds to the 2-ball of $Q(X)$. But $-2i\phi$ is the unique differential in the 2-ball of $Q(X)$ which pairs to 1 against $\frac{i}{2}\mu_\phi$.
Thus, $dF_X$ is integration against $-2i\phi$, as claimed.

By construction, the tangent vector
$\pd{\E^{\phi}}{\lambda_j}(i,\ldots, i)$ is represented by the Beltrami differential which is equal to $\frac{i}{2}\mu_{\phi}$ on the $j$th cylinder $\Pi_j$
and zero elsewhere. To obtain \eqref{Deriv}, compute
\begin{align*}
\pd{f}{\lambda_j}(i,\ldots, i) &= dF_X\lp \pd{\E^{\phi}}{\lambda_j}(i,\ldots, i) \rp\\ 
                                           &= \int_{\Pi_j} \lp -2i\phi\rp \lp \frac{i}{2} \mu_{\phi} \rp\\
                                           &= \int_{\Pi_j} \lpi \phi \rpi\\ &= a_j.
\end{align*}
\qed

The key tool in the proof of Theorem \ref{criterion} is a complex-analytic result concerning
the space $\mathcal{D}$ of holomorphic functions $\H^k \ra \H$ which satisfy condition $\eqref{Diag}$.
%Let $\mathcal{D}$ be the space of holomorphic functions $\H^k \ra \H$ satisfying condition $\eqref{Diag}$.
Consider the conjugation of action of $\R$ on $\mathcal{D}$:
$$f_t(\lambda_1,\ldots, \lambda_k) = f(\lambda_1-t,\ldots, \lambda_k-t)+t.$$
We call the map $(t,f)\mapsto f_t$ the {\em translation flow} on $\mathcal{D}$.
The translation flow is well-behaved; for any $f\in \mathcal{D}$, the orbit $\{f_t\}$ spends most of its time
close to a linear function. More precisely, we have the following result, proven in \cite{Ge}.  
 
\begin{Theorem}\label{flow}
Let $f \in \mathcal{D}$. Define $\g\in \mathcal{D}$ by
$$\mathbf{g}(\lambda_1,\ldots, \lambda_n) = \sum_{j=1}^k a_j\lambda_j,\text{where }
a_j = \pd{f}{\lambda_j}(i,\ldots, i).$$
Let $U$ be any neighborhood of $\g$ in the compact-open topology.
Then the set 
$$S=\{t\in \R | f_t\in U \}$$
has density 1 in $\R$:
$$\lim_{r\ra\infty} \frac{m(S\cap [-r,r])}{2r}=1,$$
where $m$ is the Lebesgue measure.
\end{Theorem}

We now prove the main result of this section.

{\em Proof of Theorem \ref{criterion}:}
Let $F:\T \ra \H$ be a holomorphic map such that $F\circ \tau^{\phi} = \id_\H$. Then $f = F\circ \E^\phi$ satisfies equations \eqref{Diag} and \eqref{Deriv}.
To apply Theorem \ref{flow} in the present context, we need to approximate translations in the polyplane $\E^\phi(\H^k)$ by Dehn multi-twists.
We will find a sequence $t^1,t^2,\ldots$ of real numbers
and a sequence $\alpha^1,\alpha^2,\ldots$ of mapping classes so that $t^N + F\circ \alpha^N$ converges to the desired map $G$. 

Let $m_j$ denote the modulus of the $j$th cylinder of $\phi$. Fix $\vn>0$. We claim we can choose $t\in \R$ so that
\begin{enumerate}[label = (\roman*)]
\item
The distance from $t$ to the nearest point in $m_j^{-1}\Z$ is less than $\vn$, for each $j=1,\ldots, k$.
\item
$$d\lp f_t, \g \rp < \vn,$$
where $\g(\llambda) = \sum a_j\lambda_j$ and $d$ is a fixed metric inducing the compact-open topology on $\mathcal{O}(\H^k)$.
\end{enumerate}
To see this, let $S_1$ be the set of $t\in \R$ satisfying Condition (i), and let $S_2$ be the set satisfying the Condition (ii).
By standard results on equidistribution of linear flows on the $k$-torus, the set $S_1$ has positive density in $\R$.
By Theorem \ref{flow}, the set $S_2$ has density 1. Therefore, the intersection $S_1\cap S_2$ is nonempty,
which is what we need.

So pick $t$ satisfying the above conditions and find integers $N_j$ so that $\lpi t-\frac{N_j}{m_j}\rpi < \vn$.
Let $\alpha$ be the multi-twist which twists $N_j$ times about the $j$th cylinder. Now,  set 
$$G_\vn = t + F\circ \alpha^{-1}.$$

Then 
$$d(G_\vn \circ \E^{\phi}, \g)\leq d(G_\vn\circ \E^\phi, f_t) + d(f_t,\g).$$
By Condition (ii), the second term is less than $\vn$. By Condition (i), the first term is small;
to see this, write
\begin{align*}
d_\H\lp G_\vn\circ \E^\phi(\llambda), f_t(\llambda)  \rp &= d_\H\lp f\lp \lambda_1- \frac{N_1}{m_1},\ldots, \lambda_k- \frac{N_k}{m_k} \rp, f\lp \lambda_1-t,\ldots, \lambda_k-t\rp \rp\\
&\leq \max d_\H\lp \lambda_j-\frac{N_j}{m_j}, \lambda_j-t\rp\\
&\leq \max d_\H\lp \lambda_j, \lambda_j+\vn\rp,
\end{align*}
where, in the equality, we have used \eqref{Twist} and, in the first inequality, we have used the fact that
the Kobayashi distance on $\H^k$ is the max of the Poincar\'e distances on the factors.
The last displayed quantity goes to 0 locally uniformly as $\vn \ra 0$.

Thus, the sequence $$G_{\frac{1}{2}}, G_{\frac{1}{3}}, G_{\frac{1}{4}},\ldots$$
is a normal family and any subsequential limit $G$ satisfies $d\lp G\circ \E^\phi, \g \rp = 0$, i.e. $G\circ \E^\phi= \g$ . \qed 

If $\psi = \mmu\cdot \phi$ is another element of the orbit $\H^k\phi$, then 
$\tau^{\psi}(\lambda)$ is the surface obtained by applying to the $j$th cylinder of $\phi$ the linear transformation
$$
\lp \begin{array}{cc} 1 & \Re\lp \lambda\rp \\ 0 & \Im \lp \lambda\rp \end{array}\rp \lp \begin{array}{cc} 1 & \Re \lp \mu_j\rp \\ 0 & \Im \lp\mu_j\rp \end{array}\rp =
\lp \begin{array}{cc} 1 & \Re \lbrac  \Im \lp \mu_j\rp \lambda + \Re \lp\mu_j\rp \rbrac \\ 0 & \Im \lbrac \Im \lp\mu_j\rp \lambda + \Re\lp\mu_j\rp \rbrac \end{array}\rp, 
$$
so 
$$G \circ \tau^\psi(\lambda) = \sum_j a_j\lbrac \Im \lp\mu_j\rp \lambda + \Re\lp\mu_j\rp \rbrac = c\lambda + d,$$
where $c = \sum_j \Im\lp \mu_j\rp$ and $d = \sum_j \Re \lp\mu_j\rp$. So $\frac{G-d}{c}$ is a retraction onto the disk generated by
$\psi$. Thus, we have proved
\begin{Corollary}\label{Cor}
Let $\phi$ be a Jenkin-Strebel differential with $k$ cylinders.
If $\phi$ generates a holomorphic retract, then so does every differential in its $\H^k$ orbit.  
\end{Corollary}

\section{The L-shaped pillowcase}\label{NEQ}
\begin{figure}[h]
\includegraphics[scale=1.3]{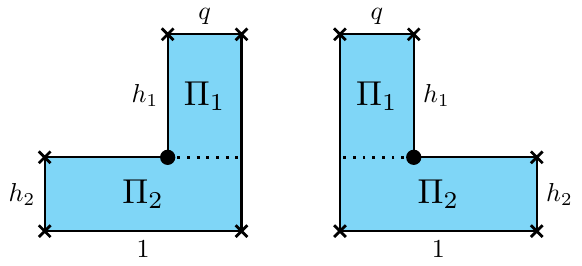}
\caption{Gluing two copies of the $L$ along corresponding edges yields the differential $\phi(h_1,h_2,q)$ on $S_{0,5}$. The crosses indicate
poles and the dot indicates a simple zero.}
\label{Lshaped}
\end{figure}

\begin{figure}[h]
\includegraphics[scale=1.3]{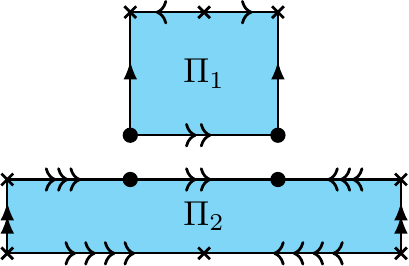}
\caption{The cylinders of the $L$-shaped pillowcase.}
\label{Lcylinder}
\end{figure}

Doubling a right-angled $L$-shaped hexagon along its boundary yields a Jenkins-Strebel differential on $S_{0,5}$
called an {\em L-shaped pillowcase} (Figures \ref{Lshaped}, \ref{Lcylinder}). In this section, we sketch the proof \cite{Mar} that an $L$-shaped pillowcase does not generate a holomorphic retract, and thus that the Carath\'eodory and Kobayashi metrics on $\T_{0,5}$ do not coincide. We then prove Theorem \ref{neq2}, which states the two metrics on $\T_{g,n}$ are different whenever $\dim_\C\T_{g,n}\geq 2$.

The $L$-shaped pillowcase has two cylinders, $\Pi_1$ and $\Pi_2$.
We let $h_i$ denote the heights of $\Pi_i$. We denote the length of $\Pi_1$ by $q$ and normalize so that the length
of $\Pi_2$ is 1. We call the resulting quadratic differential $\phi(h_1, h_2, q)$ and its underlying marked surface $X(h_1,h_2, q)$.

\begin{Theorem}\label{Not}\cite{Mar}
The Teichm\"uller disk generated by $\phi(h_1, h_2, q)$ is not a holomorphic retract.
\end{Theorem}
{\em Sketch:}
By Corollary \ref{Cor}, it suffices to show that $\phi_0 := \phi(1,1,q)$ does not generate a retract. Suppose
to the contrary that the disk $\tau^{\phi_0}$ is a holomorphic retract of $\T_{0,5}$. Then by Theorem
\ref{criterion}, there is a holomorphic $G:\T_{0,5} \ra \H$ so that
$$G\circ \E^{\phi_0} = a_1\lambda_1 + a_2\lambda_2,$$
with $a_1= \frac{q}{1+q}$ and $a_2 = \frac{1}{1+q}$.

\begin{figure}[h!]
\centering
\subfloat[][$X(0, h_2, q)$]{\includegraphics[scale=2]{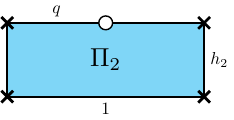}\label{CollapsedL}}
\hfil
\subfloat[][$\gamma(t)\mkern-1.3mu = \mkern-1.3mu X(0, h_2, q-t)$]{\includegraphics[scale=2]{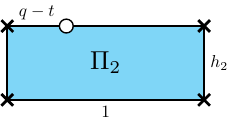}\label{CollapsedLpath}}
\caption{Collapsing the top cylinder of an L-shaped pillowcase yields a well-defined quadratic differential. A simple pole and zero ``cancel" to give
a regular point at a puncture. We get a path in $\T_{0,5}$ by moving this puncture horizontally.} 
\end{figure}

The idea is to reach a contradiction by examining the regularity of $G$ at the boundary of $\mathcal{E}^{\phi_0}(\H^2)$.
To this end, note that the differential $\phi(0, 1, q)$ obtained by collapsing $\Pi_1$ is a well-defined element of $\widetilde{\mathcal{Q}}_{0,5}$ 
(Figure \ref{CollapsedL}).
Now observe that $\gamma(t)= X(0,1, q-t)$ is a smooth path in Teichm\"uller space (Figure \ref{CollapsedLpath}). 
Since $G$ is holomorphic, $G\circ \gamma$ is a smooth path in $\H$.

Now, an argument using the Schwarz-Christoffel mappings shows that $X(0,1,q-t)$ is in $\mathcal{E}^\phi(\H^2)$ for each
$t\in (0,q)$. In fact, there is a unique pair of positive numbers $h_1(t), h_2(t)$ so that  $$X(0,1,q-t) = X(h_1(t), h_2(t), q ).$$    

Thus,
\begin{align*}
G\circ \gamma(t) &= G\circ E^\phi (h_1(t)i, h_2(t)i)\\
                              &= \lbrac a_1h_1(t) + a_2h_2(t)\rbrac i.
\end{align*}
An involved computation with the Schwarz-Christoffel mappings determines the asymptotics of $h_1(t), h_2(t)$ for small positive $t$.
See \cite{Mar} Sections 8,9 for details.
The end result is that there are constants $\beta_1$ and $\beta_2\neq 0$ so that
$$G\circ \gamma(t) -  G\circ \gamma(0) =  \beta_1(1+o(1)) \frac{t}{\log t^{-1}} + \beta_2(1+o(1))\frac{t^2}{\log t^{-1}}  + o\lp\frac{t^2}{\log t^{-1}}\rp,$$
which is incompatible with the fact that $G\circ \gamma$ is thrice-differentiable. \qed

{\em Proof of Theorem \ref{neq2}:}
We have seen that the Kobayashi and Carath\'eodory metrics on $\T_{0,5}$ are different. 
To prove the corresponding fact for $\T_{g,n}$, we note that there is an embedding $\iota:\T_{0,5} \ra \T_{g,n}$
which is holomorphic and isometric for the Kobayashi metric. 
The embedding is constructed using an elementary covering argument; we give the details in
the Appendix.

Now, let $\phi\in \widetilde{\mathcal{Q}}_{0,5}$ be an $L$-shaped pillowcase.
Then $\iota\circ \tau^{\phi}$ is a Teichm\"uller disk
in $\T_{g,n}$. A holomorphic retraction $F:\T_{g,n} \ra \H$ onto this disk would yield a retraction
$F\circ \iota$ onto $\tau^{\phi}$, contradicting Theorem \ref{Not}. So $\iota\circ \tau^{\phi}$ is not a retract
and thus the two metrics on $\T_{g,n}$ are different.
\qed

\newpage
\section{The five-times punctured sphere and twice-punctured torus}\label{MAIN}
\subsection{Jenkins-Strebel differentials on $S_{0,5}$}
In this section, we classify the Jenkins-Strebel differentials on the five-times punctured sphere $S_{0,5}$.
We establish the following:
\begin{Proposition}\label{Fivehole}
~\newline
\vspace*{-.8cm}
\begin{enumerate}[label=\normalfont(\alph*), ref=(\alph*)]
\item
Any differential with an odd-order zero has a Jenkins-Strebel differential with two cylinders
in its $\GL_2^+$ orbit closure.
\item 
Any differential with two cylinders has an $L$-shaped pillowcase in its $\H^2$-orbit.
\end{enumerate}
\end{Proposition}
\vspace{-10pt}
Combined with Proposition \ref{Closure}, Corollary \ref{Cor}, and Theorem \ref{Not}, this proves our main result that Conjecture \ref{Conj} holds
in the case of $\T_{0,5}$. We summarize the proof of the main result in Section \ref{M}.
   
An integrable holomorphic quadratic differential $\phi$ on $S_{0,5}$ either has five poles and a simple zero
or four poles and no zeros. In the second case, $\tau^{\phi}$ is a retract by Theorem \ref{Kra}. (Alternatively, note that the forgetful map $F:\T_{0,5} \ra \T_{0,4}$ which ``fills in" the puncture at the regular point restricts to a biholomorphism $F \circ \tau^{\phi}:\H \ra \T_{0,4}$.) 

\begin{figure}[h!]
\centering
\subfloat[][]{\includegraphics[scale=.9]{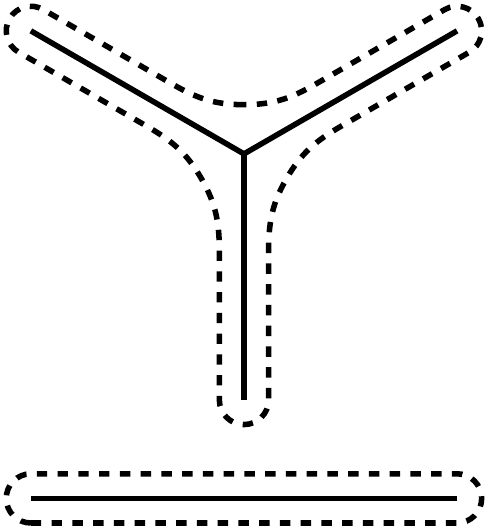}\label{CriticalOne}}
\hfil
\subfloat[][]{\includegraphics[scale=.9]{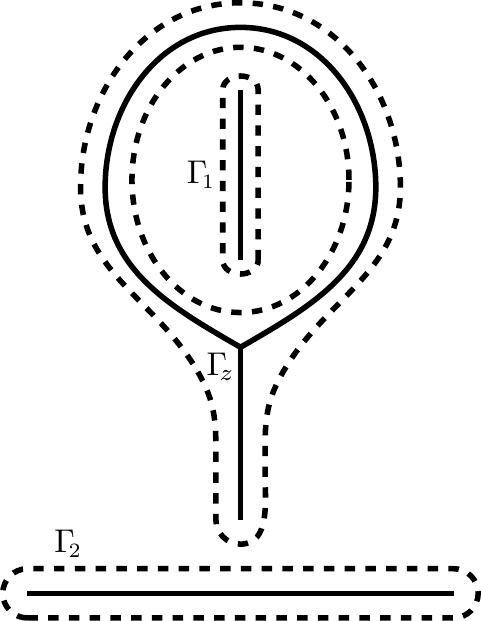}\label{CriticalTwo}}
\vspace{-5pt}
\caption{If $\phi\in \mathcal{Q}_{0,5}$ has a simple zero, then its critical graph is of one of the two indicated types.} 
\end{figure}

Conjecture \ref{Conj} for $\T_{0,5}$ reduces to showing that, in case $\phi$ has a simple zero, $\tau^{\phi}$ is not a retract.
As we observed in Section \ref{Horocycle}, it suffices to consider the case that $\phi$ is Jenkins-Strebel.
Let $\Gamma$ be the critical graph of $\phi$, i.e. the union of the horizontal rays emanating from zeros and poles.
Since $\phi$ is Jenkins-Strebel, each horizontal ray emanating from a singularity ends at a singularity.
Thus, $\Gamma$ is a finite graph on $S^2$ with a valence three vertex at the simple zero
and valence one vertices at each of the poles. (In this context, the valence of a vertex is the number of {\em half-edges} incident on it; a loop counts twice
towards the valence of the incident vertex.) Each boundary component of an $\vn$-neighborhood of $\Gamma$ is a closed horizontal curve.
There are two possibilities:

{\bf (Case 1) All three of the horizontal rays emanating from the zero terminate at simple poles.} (See Figure \ref{CriticalOne}.)
In this case, the two remaining poles are joined by a horizontal segment. The $\vn$-neighborhood
of $\Gamma$ has two boundary components. Thus, $\phi$ is Jenkins-Strebel with one cylinder.

\begin{figure}[h!]
\includegraphics[scale=.8]{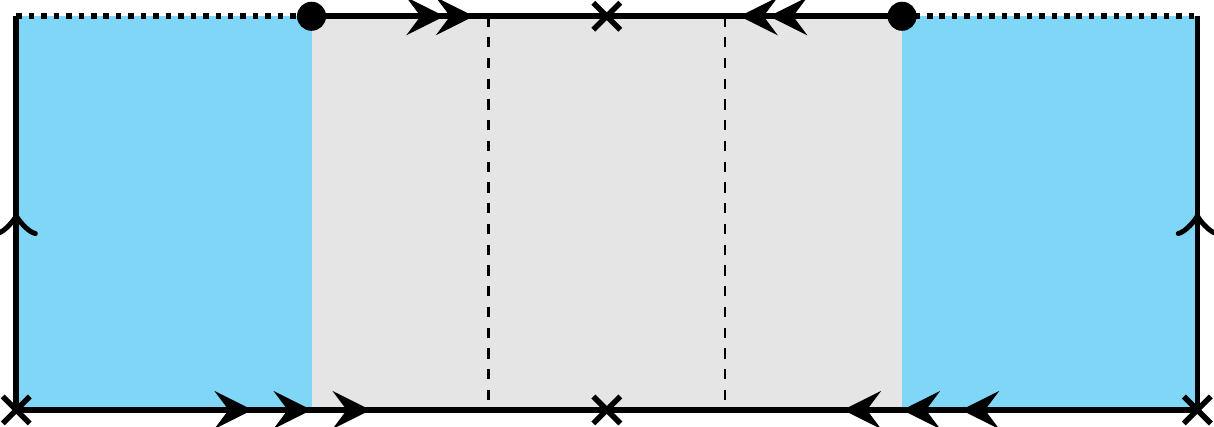}
\caption{A Strebel differential on $S_{0,5}$ with one cylinder and an odd-order zero. Crosses indicate poles,
and the large dot indicates a zero. The dashed curve is the core of a closed vertical cylinder.}
\label{StrebelOne}
\end{figure}

By shearing $\phi$ appropriately, we may assume that there is a vertical geodesic connecting a pole
in one component of $\Gamma$ to a pole in the other component.
Then $\phi$ has a vertical cylinder which is not dense in $S_{0,5}$ (Figure \ref{StrebelOne}).
So by Theorem \ref{hclosure}(b), the horocycle orbit closure of the rotated differential 
$$\lp \begin{array}{cc} 0 & -1\\ 1 &0\end{array} \rp \cdot \phi$$
contains a Jenkins-Strebel differential with two cylinders.

{\bf (Case 2) One of the horizontal rays emanating from the zero terminates at the zero.} (See Figure \ref{CriticalTwo}.)
Another ray emanating from the zero terminates at a pole.
Let $\Gamma_z$ denote the component of $\Gamma$ containing the zero.
Let $M_z$ denote the $\vn$-neighborhood of $\Gamma_z$. 
There are two remaining pairs of poles; each pair is joined by a horizontal segment.
Let $\Gamma_1, \Gamma_2$ denote these horizontal segments. Let $M_1, M_2$
denote their $\vn$-neighborhoods.

The boundary $\partial M_z$ has two components, one longer than the other.
Thus, one of the components of $\partial M_z$, say the shorter one, is homotopic to $\partial M_1$, and the other is homotopic to $\partial M_2$. 
So $\phi$ is Jenkins-Strebel with two cylinders.

Now, shear a cylinder of $\phi$ so that there is a vertical segment connecting the zero in $\Gamma_z$ and a pole in 
$\Gamma_1$. Next, shear the other cylinder so that there is a vertical geodesic
connecting the pole of $\Gamma_z$ to a pole of $\Gamma_2$. The resulting differential is an $L$-shaped pillowcase in the same
$\H^2$-orbit as $\phi$.

{\em Proof of Proposition \ref{Fivehole}:}
To prove the first part, suppose $\phi \in \mathcal{Q}_{0,5}$ has an odd-order zero. Let $\psi\in \ol{\GL_2^+\phi}$ be a Jenkins-Strebel differential with a simple zero.
If $\psi$ has two cylinders, we are done. If it has only one, then we are in {\bf Case 1}, so $\psi$ has a differential
with two cylinders in its orbit closure.

To prove the second part, suppose $\phi \in \widetilde{\mathcal{Q}}_{0,5}$ has two cylinders. Then we are in {\bf Case 2}, so $\phi$ is in the same $\H^2$
orbit as an $L$-shaped pillowcase. \qed

\captionsetup[subfigure]{position=t,singlelinecheck=off}
 
\subsection{Proof of the main theorem}\label{M}
We collect our results and classify Carath\'eodory geodesics in $\T_{0,5}$ and $\T_{1,2}$.

{\em Proof of Theorem \ref{Main}:}
We first prove the result for $\T_{0,5}$.
The ``if" direction follows from Kra's Theorem \ref{Kra}. For the ``only if" direction, let $\phi \in {\mathcal{Q}}_{0,5}$
be a differential with an odd-order zero. Suppose for the sake of contradiction that $\phi$ generates a holomorphic retract.
By Theorem \ref{hclosure} and Proposition \ref{Fivehole}(a), the $\GL_2^+$ orbit closure of $\phi$ contains a Jenkins-Strebel differential $\phi'$
with two cylinders and an odd-order zero. By Proposition \ref{Fivehole}(b), the orbit $\H^2\cdot \phi'$ contains an $L$-shaped pillowcase $\phi''$.
By Proposition \ref{Closure} and Corollary \ref{Cor}, $\tau^{\phi''}$ is a holomorphic retract. This contradicts Theorem \ref{Not}. 

To prove the result for $\T_{1,2}$, recall that there is an orientation-preserving involution 
$\alpha \in \MCG_{1,2}$ which fixes every point of $\T_{1,2}$.
For each $X \in \T_{1,2}$, the class $\alpha$ is represented by a conformal involution of $X$
which fixes four points and swaps the punctures. The quotient of $X$ by the involution
is a surface $X'$ of genus zero with five marked points. 
There is a unique complex structure on $X'$ making the quotient $f:X\ra X'$ 
a holomorphic double cover branched over four marked points.

The map $\T_{1,2} \ra \T_{0,5}$ sending $X$ to $X'$ is a biholomorphism.
Let $\phi$ be a quadratic differential on $X'$, and let $f^*\phi$ be its pullback to $X$.
Then the Teichm\"uller disk $\tau^{\phi}$ in $\T_{0,5}$ corresponds to the
disk $\tau^{f^*\phi}$ in $\T_{1,2}$.
To complete the proof, we observe that $f^*\phi$ has an odd-order zero if and only if $\phi$ does.
Indeed, a simple zero of $\phi$ is necessarily unramified and thus lifts to two simple zeroes of $f^*\phi$.
On the other hand, ramified poles and ramified regular points lift to regular points and double zeroes, respectively.
So if $\phi$ has no odd-order zeros, neither does $f^*\phi$.
 \qed    

\section{Appendix: $\T_{0,5}$ embeds in $\T_{g,n}$}
\begin{figure}[h]
\centering
\subfloat[][]{\includegraphics[scale=2.5]{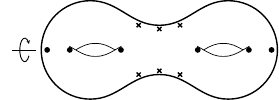}\label{Involution}}
\quad
\subfloat[][]{\includegraphics[scale=2.5]{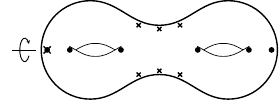}\label{InvolutionOdd}}
\caption{The desired involution is 180 degree rotation about the indicated axis. Crosses indicate marked points. Dots indicate fixed points of the involution.}
\end{figure}

\vspace{-13pt}  
We used the following Lemma in Section \ref{NEQ} to deduce that the Kobayashi and Carath\'eodory metrics on $\T_{g,n}$
are different.
\begin{Lemma}
Whenever $\dim_\C \T_{g,n} = 3g-3+n \geq 2$, there is a holomorphic and isometric embedding
$\T_{0,5} \ra \T_{g,n}$.
\end{Lemma}
{\em Proof:}
If $n=2m$ is even, then $S_{g,n}$ admits an involution $\alpha$ fixing $2g+2$ points, none of which are marked (Figure \ref{Involution}).
If $n=2m+1$ is odd, then $S_{g,n}$ admits an involution $\alpha$ fixing $2g+2$ points, one of which is marked (Figure \ref{InvolutionOdd}).
In either case, the quotient is $S_{0, n_1}$ where $n_1 = m + 2g+2\geq 5$.
The space $\T_{0,n_1}$ embeds holomorphically and isometrically in $\T_{g,n}$ as the fixed point set of the action of $\alpha$ 
(see e.g. \cite{Fa} p. 370). (In case $g=1$ and $n=2$, we obtain the isomorphism $\T_{0,5} \cong \T_{1,2}$ described in the last section.
This construction also yields the other two coincidences $\T_{0,4}\cong \T_{1,1}$ and $\T_{0,6} \cong \T_{2,0}$.)
If $n_1 = 5$, we are done. Otherwise, by the same construction as above, the surface $S_{0,n_1}$ admits an involution
with quotient $S_{0, n_2}$ and $5\leq n_2 < n_1$ (Figure \ref{InvolutionSphere}). So we have embeddings $$\T_{0,n_2}\hra \T_{0,n_1} \hra \T_{g,n}.$$
Continuing inductively, we get an embedding $\T_{0,5} \hra \T_{g,n}$.   
\qed

%\floatsetup[figure]{style=plain,subcapbesideposition=top}
\setlength{\belowcaptionskip}{-10pt}
\begin{figure}[h]
\vspace{-20pt} 
\centering
\subfloat[][]{\includegraphics[scale=2.5]{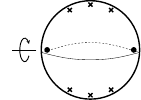}\label{InvolutionSphere}}
\hspace{110pt}
\subfloat[][]{\includegraphics[scale=2.5]{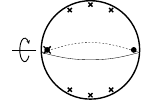}\label{InvolutionOdd}}
\caption{}
\label{InvolutionSphere}
\end{figure}

\pagebreak
\section{Acknowledgements}
We would like to thank Alex Wright for helpful discussion.


\begin{thebibliography}{}

\bibitem{An1}
S. Antonakoudis, \textit{Isometric disks are holomorphic}. Invent. Math. 207(3)  (2017), 1289-1299.

\bibitem{An2}
S. Antonakoudis, \textit{Teichm\"uller spaces and bounded symmetric domains do not mix isometrically}, GAFA 27 (2017), 453-465.

\bibitem{Ea}
C. Earle, \textit{On the Carath\'eodory metric in Teichm\"uller spaces.} Discontinuous groups and
Riemann surfaces, 99-103. Ann. of Math. Studies, No. 79, Princeton Univ. Press, Princeton, N.J. (1974)

\bibitem{Fa}
B. Farb, and D. Margalit. \textit{A Primer on Mapping Class Groups.} Princeton University Press, 2011.

\bibitem{Ge}
D. Gekhtman, \textit{Asymptotics of the translation flow on holomorphic maps out of the polyplane.} arXiv:1702.02177 (2017).

\bibitem{Gu}
S. Gupta and H. Seshadri, \textit{On domains biholomorphic to Teichm\"{u}ller spaces.} arXiv:1701.06860 (2017).

\bibitem{Je}
J. Jenkins, \textit{On the existence of certain general extremal metrics}. Ann. of Math. 66 (1957),
440-453.

\bibitem{Ja}
M. Jarnicki and P. Pflug, \textit{Invariant Distances and Metrics in Complex Analysis}. Walter de Gruyter and Co., Berlin (1993).

\bibitem{Ke}
S. Kerckhoff. \textit{The asymptotic geometry of Teichmuller space.} Topology 19.1 (1980), 23-41.

\bibitem{Kr}
I. Kra, \textit{The Carath\'eodory metric on abelian Teichm\"uller disks}. Journal Analyse Math. 40 (1981), 129-143.

\bibitem{Le1}
L. Lempert, \textit{La m\'etrique de Kobayashi et la repr\'esentation des domaines sur la boule},
Bull. Soc. Math. France 109 (1981), 427-474.

%\bibitem{Le2}
%L. Lempert, \textit{Holomorphic retracts and intrinsic metrics in convex domains}. Analysis
%Mathematica 8 (1982), 257-261.

\bibitem{Mar}
V. Markovic, \textit{Carath\'eodory's metric on Teichm\"uller spaces and L-shaped pillowcases}. To appear in Duke Math. J.

\bibitem{Mas}
H. Masur, \textit{Interval exchange transformations and measured foliations}. Ann. of Math. 115.1 (1982), 169-200.

\bibitem{Mi}
Y. Minsky and B. Weiss, \textit{Nondivergence of horocyclic flows on moduli space}. Journal f\"ur die Reine und Angewandte Mathematik 552 (2002), 131- 177.

\bibitem{Ro}
H. Royden, \textit{Automorphisms and isometries of Teichm\"uller space}. Advances in the Theory of Riemann Surfaces (Proc. Conf., Stony Brook, N.Y., 1969) pp. 369-383 Ann. of Math. Studies, No. 66. Princeton Univ. Press, Princeton, N.J. (1971).

\bibitem{Sm}
J. Smillie and B. Weiss, \textit{Minimal sets for flows on moduli space}. Israel J. of Math. 142 (2004), 249-260.

\bibitem{Ta}
H. Tanigawa, \textit{Holomorphic mappings into Teichm\"uller spaces.} Proc. Amer. Math. Soc. 117 (1993), 71-78.

\bibitem{Ve}
W. Veech. \textit{Gauss measures for transformations on the space of interval exchange
transformations}. Ann. of Math. 115 (1982), 215-242.

\bibitem{Wr}
A. Wright. \textit{Cylinder deformations in orbit closures of translation surfaces}. Geometry \& Topology 19.1 (2015), 413-438.


\end{thebibliography}
\end{document}